\documentclass[a4paper]{article}

\usepackage{geometry,amsmath,amssymb,commath,graphicx,bm}
\usepackage{subcaption}
\usepackage{hyperref}
\usepackage{cleveref}
\usepackage[title]{appendix}
\usepackage{doi}
\usepackage[square,numbers]{natbib}
\hypersetup{pdfborder = 0 0 0, colorlinks=false, linkcolor=black,citecolor=black, filecolor=black, urlcolor=black}
\usepackage{authblk}
\usepackage{geometry}
\usepackage{bm}
\usepackage{commath}
\usepackage{graphicx}
\usepackage{todonotes}

\newtheorem{thm}{Theorem}

\newtheorem{rmk}[thm]{Remark}      

\newcommand{\Res}{\operatorname{Res}}
\newcommand{\resbrac}[1]{\Res\left[#1\right]}
\renewcommand{\Re}{\operatorname{Re}}
\renewcommand{\Im}{\operatorname{Im}}

\newcommand{\pars}[1]{\left( #1 \right)}%
\newcommand{\preimage}{t^*}

\title{Singularity swap quadrature for nearly singular line integrals
  on closed curves in two dimensions}

\author{Ludvig af Klinteberg\thanks{E-mail: \texttt{ludvig.af.klinteberg@mdu.se}}}
\affil{Department of Mathematics and Physics,\\M\"alardalen University,\\V\"aster{\aa}s, Sweden}
\date{}

\begin{document}
\maketitle

\begin{abstract}
  This paper presents a quadrature method for evaluating layer
  potentials in two dimensions close to periodic boundaries,
  discretized using the trapezoidal rule.  It is an extension of the
  method of singularity swap quadrature, which recently was introduced
  for boundaries discretized using composite Gauss-Legendre
  quadrature.  The original method builds on swapping the target
  singularity for its preimage in the complexified space of the curve
  parametrization, where the source panel is flat.  This allows the
  integral to be efficiently evaluated using an interpolatory
  quadrature with a monomial basis.
  In this extension, we use the target preimage to swap the
  singularity to a point close to the unit circle. This allows us to
  evaluate the integral using an interpolatory quadrature with complex
  exponential basis functions. This is well-conditioned, and can be
  efficiently evaluated using the fast Fourier transform.
  The resulting method has exponential convergence, and can be used to
  accurately evaluate layer potentials close to the source geometry.
  We report experimental results on a simple test geometry, and
  provide a baseline Julia implementation that can be used for further
  experimentation.
\end{abstract}

\section{Introduction}

In integral equation methods, one of the long-standing challenges is
the evaluation of layer potentials for targets points close to the
source geometry. The integrals by which these layer potentials are
computed are referred to as nearly singular, since the integrand is
evaluated close to a singularity. This near singularity reduces the
smoothness of the integrand, thereby severely reducing the accuracy of
common quadrature rules such as Gauss-Legendre and the trapezoidal
rule. To tackle this, many specialized quadrature methods have been
proposed for the case of the domain being a one-dimensional line. For
panel-based Gauss-Legendre quadrature in the plane, the high-order
method of \cite{Helsing2008} (often referred to as the Helsing-Ojala
method) is efficient and well-established. For the trapezoidal rule,
which is exponentially convergent on a closed curve, several
fixed-order methods have been proposed that reduce the impact of the
near singularity by means of a local correction, see
e.g. \cite{Beale2001, Carvalho2018, Nitsche2022}. In addition,
exponentially convergenent near evaluation is available through the
``globally compensated'' quadrature method introduced in
\cite{Helsing2008} and extended in \cite{Barnett2015}.

Recently, the method of singularity swap quadrature (SSQ) was proposed
in \cite{AfKlinteberg2020line}. The method evaluates nearly singular
line quadratures in two and three dimensions by first ``swapping'' the
singularity from the target point in the plane to the preimage of the
target point in the curve's parameterization. After the swap, the
integral can be evaluated using interpolatory quadrature with a
monomial basis, following the method of \cite{Helsing2008}. This was
shown to improve the accuracy for curved panels compared to
\cite{Helsing2008}, and allowed the method to be extended to line
integrals in three dimensions.

The SSQ method was in \cite{AfKlinteberg2020line} described for curves
discretized using panel-based Gauss-Legendre quadrature. In this brief
follow-up, we show how the method can be extended to closed curves in
the plane that have been discretized using the trapezoidal rule,
resulting in a global correction with exponential convergence. In
doing so, we assume that the reader has at least some familiarity with
the original method.

\section{Method}

Since we are considering the problem in two dimensions, let us
identify $\mathbb R^2$ with $\mathbb C$, and let $\Gamma$ be a closed
curve parameterized by a smooth function
$\gamma(t) = g_1(t) + i g_2(t)$, with $t \in [0, 2\pi)$. Furthermore,
let $\sigma$ be a smooth function defined on $\Gamma$, which we will
refer to as the \emph{density}. For a point $z\in\mathbb C$, typically
lying close to $\Gamma$, our integrals of interest are here
\begin{align}
  I_L &= I_L(z) := \oint_\Gamma \sigma \log\abs{\tau-z}\abs{\dif\tau},
        \label{eq:complex_L}\\
  I_m &= I_m(z) :=
        \oint_\Gamma \frac{\sigma \dif\tau}{(\tau-z)^m} 
        , 
        \qquad m=1,2,\dots .
        \label{eq:complex_m}
\end{align}
Being periodic, these integrals are suitable for discretization using
an N-point trapezoidal rule, due to its exponential convergence
\cite{Trefethen2014}. The quadrature nodes are then $\gamma(t_j)$,
$j=0,\dots,N-1$, with $t_j = 2\pi j/N$ and the quadrature weights
$w_j=2\pi/N$. We will now outline how these discretized integrals can
be evaluated using SSQ. Throughout, we will assume that we only have
access to the discrete values of the functions $\sigma$, $\gamma$ and
$\gamma'$ at the nodes, and not to their analytic definitions.

\subsection{The Cauchy integral}
We begin with the simplest (and perhaps most common) case, the Cauchy integral,
\begin{align}
  I_1(z) &= \oint_\Gamma \frac{\sigma \dif \tau}{\tau - z}
           = \int_0^{2\pi} \frac{\sigma(t)\gamma'(t)\dif t}{\gamma(t)-z} .
           \label{eq:cauchy_int}
\end{align}
Let us now assume that $z$ is close to $\Gamma$, which corresponds to
the integrand in \eqref{eq:cauchy_int} having a singularity close to
the real line at the preimage $\preimage \in \mathbb C$ such that
\begin{align}
  \gamma(\preimage) = z.
  \label{eq:preimage_equation}
\end{align}
In \cite{Barnett2014} and \cite{afKlinteberg2022}, it is shown that this
significantly reduces the convergence rate of the trapezoidal rule,
with the error being proportional to $e^{-N |\Im{\preimage}|}$. In
order to reduce the error for small $\Im \preimage$, we will now show
how to evaluate $I_1(z)$ using SSQ.

Step one is to find the preimage $\preimage$, using Newton's method
applied to a Fourier expansion of $\gamma$. We apply a Fast Fourier
Transform (FFT) to the node values $\gamma(t_0),\dots,\gamma(t_{N-1})$
to get the truncated Fourier series
\begin{align}
  \gamma(t) &\approx \sum_{k=-K}^{K} \hat\gamma_k e^{ikt}, \\
  \gamma'(t) &\approx \sum_{k=-K}^K ik\hat\gamma_k e^{ikt},
\end{align}
where $K=\frac{N-1}{2}$, assuming $N$ odd, and
$\hat\gamma_k=\mathcal F[\gamma](k)$ are the Fourier coefficients of
$\gamma$. This has a one-time cost of $\mathcal O(N\log N)$, for a
given discretization.  Substituting this expansion into
\eqref{eq:preimage_equation}, we can find $\preimage$ at an
$\mathcal O(N)$ cost using Newton's method (assuming $\mathcal O(1)$
iterations to convergence). A good initial guess for $\preimage$ is
the corresponding $t$-value of the grid point on $\Gamma$ that is
closest to $z$,
\begin{align}
  \preimage \approx \underset{t_j}{\operatorname{arg\,min}} \abs{\gamma(t_j) - z}
  .
\end{align}
We assume that $\Gamma$ is parametrized in a counter-clockwise
direction, so that $\Im \preimage > 0$ implies $z$ interior to $\Gamma$, and
vice versa.

Step two is to swap out the singularity using a function that
regularizes the integrand, and leads to an integrable
singularity. Contrary to the Gauss-Legendre panel case, we can not
regularize the denominator $\left(\gamma(t)-\gamma(\preimage)\right)^{-1}$
using the factor $(t-\preimage)$, as the latter is not periodic. Instead, we
use $(e^{it}-e^{i\preimage})$, which corresponds to swapping the singularity
to a point close to the unit circle (this will be evident
shortly). Then,
\begin{align}
  I_1(z) &= 
           \int_0^{2\pi} 
           \underbrace{
           \frac{\sigma(t)\gamma'(t) \pars{e^{it}-e^{i\preimage}}}
           {\gamma(t)-\gamma(\preimage)} 
           }_{f(t,\preimage)}
           \frac{ \dif t}{e^{it}-e^{i\preimage}}
           .  
\end{align}
The function $f$ is now the original integrand with the singularity at
$\preimage$ removed, and we expect it to be smooth. Since its values are
known at the equidistant nodes $t_0, \dots, t_{N-1}$, we use an FFT to
compute its truncated Fourier series (at an $\mathcal O(N\log N)$
cost), such that
\begin{align}
  I_1(z) &\approx
           \sum_{k=-K}^K
            \hat f_k(\preimage)
            \underbrace{
            \int_0^{2\pi} 
            \frac{ e^{ikt} \dif t}{e^{it}-e^{i\preimage}}            
            }_{p_k(\preimage)}
            = \bm{\hat f}\cdot\bm{p}
           .
           \label{eq:I1_reg}
\end{align}
This can be interpreted as our integral being approximated by an
interpolatory quadrature on the unit circle, which we can evaluate to
high accuracy as long as the coefficients $\hat f_k$ decay
sufficiently fast. The integrals $p_k$ can be evaluated analytically
(at an $\mathcal O(N)$ cost)
by rewriting them as contour integrals on the unit circle $S$. Let
\begin{align}
  \xi(t) &= e^{it},  \\
  \xi'(t) &= ie^{it},\\
  \zeta &=  e^{i\preimage} .          
\end{align}
Then $e^{ikt} = \xi(t)^k$ and
\begin{align}
  p_k(\preimage) &= 
             \frac{1}{i}
             \int_0^{2\pi} 
             \frac{ e^{i(k-1)t} i e^{it} \dif t}{e^{it}-e^{i\preimage}} \\
           &=
             \frac{1}{i}
             \int_0^{2\pi}
             \frac{\xi(t)^{k-1} \xi'(t) \dif t}
             {\xi(t) - \zeta} \\
           &= 
             \frac{1}{i}
             \oint_S \frac{\xi^{k-1} \dif\xi}{\xi-\zeta} .
\end{align}
We can evaluate this integral using the residue theorem. If
$|\zeta|<1$ (or equivalently $\Im \preimage>0$, corresponding to $z$ being an
interior point), then the integrand has a simple pole at $\zeta$, with
residue
\begin{align}
  \resbrac{\frac{\xi^{k-1}}{\xi-\zeta}, \zeta}
  = \zeta^{k-1}, \quad \mbox{ if } |\zeta|<1.
\end{align}
In addition, if $k \le 0$, then the integrand has a pole of order
$1-k$ at the origin, with residue
\begin{align}
  \resbrac{\frac{\xi^{k-1}}{\xi-\zeta}, 0}
  = -\zeta^{k-1}, \quad \mbox{ if } k \le 0.
\end{align}
The two residues cancel for $k \le 0$ and $|\zeta|<1$ (i.e.
$\Im \preimage>0$), so we get
\begin{align}
  p_k(\preimage) = 
  2\pi e^{i(k-1)\preimage}
  \begin{cases}
     1, & \text{if } \Im \preimage>0 \text{ and } k \ge 1,\\
    -1, & \text{if } \Im \preimage<0 \text{ and } k \le 0,\\
    0, & \text{otherwise}.
  \end{cases}
\end{align}
This completes the method, which has a total cost of
$\mathcal O(N + N\log N)$ per target point.

\begin{figure}[h]
  \centering
  \begin{subfigure}{.49\textwidth}
    \includegraphics[width=\textwidth]{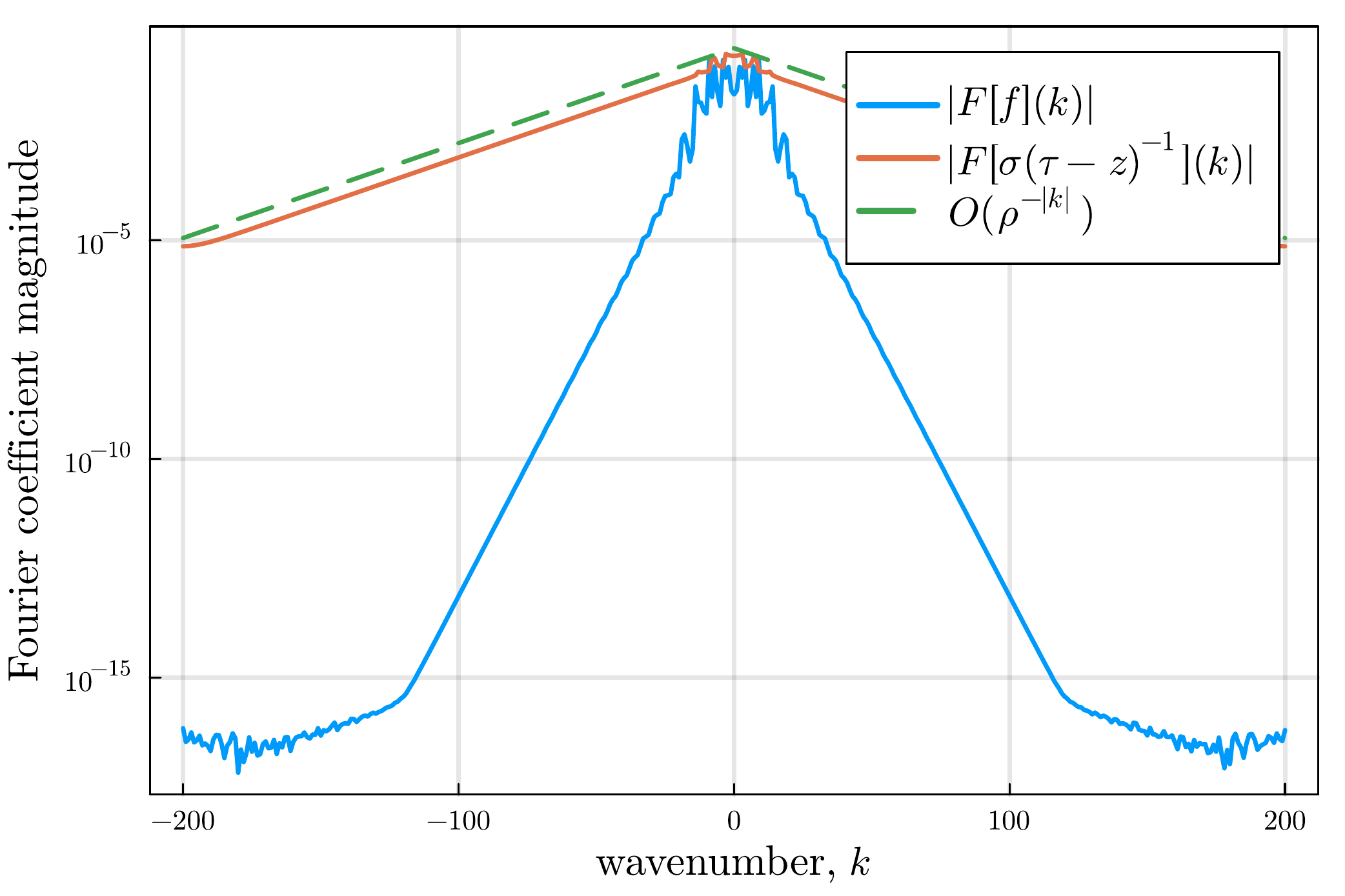}
    \caption{$\preimage=1+0.05i$}
    \label{fig:decay1}
  \end{subfigure}
  \hfill
  \begin{subfigure}{.49\textwidth}
    \includegraphics[width=\textwidth]{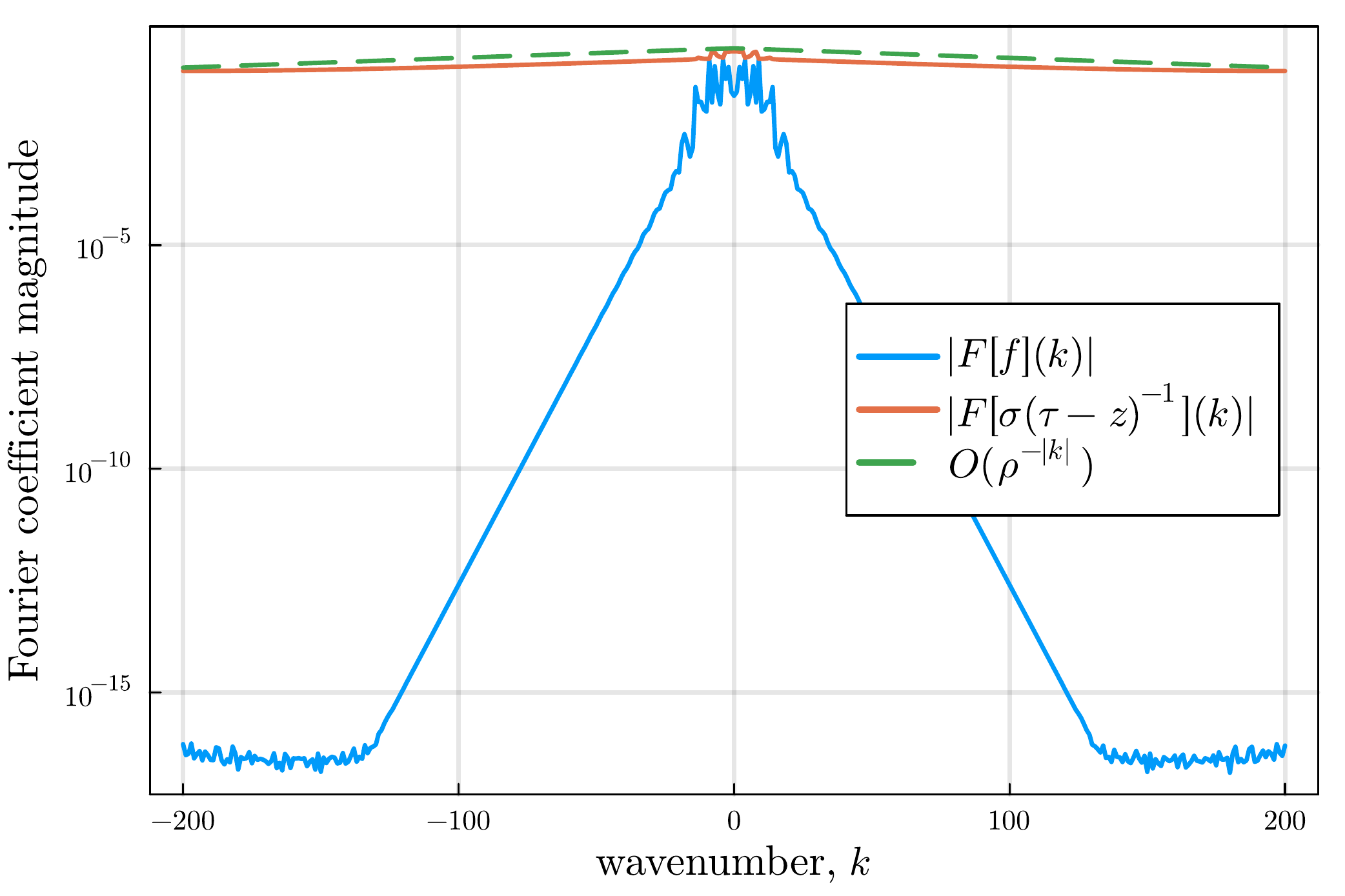}
    \caption{$\preimage=1+0.005i$}
    \label{fig:decay2}
  \end{subfigure}
  \caption{Example results for the Cauchy integral on the starfish
    geometry and density from \cref{sec:convergence}, for two
    different $\preimage$. The Fourier coefficients of the integrand
    in \eqref{eq:cauchy_int} decay approximately as $\rho^{-|k|}$,
    where $\rho=e^{|\Im t_0|}$, which gets very slow as
    $z\to\Gamma$. On the other hand, the coefficients of the
    regularized function $f$ in \eqref{eq:I1_reg} decay with a rate
    that changes little as $z\to\Gamma$.  }
  \label{fig:decay}
\end{figure}

The method can perhaps be understood in terms of Fourier
coefficents. The accuracy of the trapezoidal rule depends on the
regularity of the integrand, which in term is reflected in the decay
of the Fourier coefficients of the integrand. As seen in
\cref{fig:decay}, these appear to decay as $e^{-|k\Im\preimage|}$. In
constrast, the coefficients $\hat f_k(\preimage)$ decay rapidly and
with a rate that stays nearly constant as $z\to\Gamma$.

A weakness of the method is that in order to find $\preimage$, one
must evaluate the analytic continuation of $\gamma$ using a truncated
Fourier series, which naturally is most accurate close to
$\Gamma$. For $z$ far away from $\Gamma$ the rootfinding process will
struggle, and the accuracy of SSQ can deteriorate. However, this
happens at ranges where the base trapezoidal rule is sufficiently
accurate.

\begin{rmk}
  \label[remark]{remark:direct}
  The reader may note that one could in fact apply a variant of the
  Helsing-Ojala quadrature directly to (\ref{eq:cauchy_int}).
  First solve interpolation problem
  \begin{align}
    \sum_{k=-K}^K c_k \tau_i^k = \sigma_i, \quad i=0,\dots,N-1,
  \end{align}
  where $\tau_i=\gamma(t_i)$ and $\sigma_i=\gamma(t_i)$. Then evaluate
  the layer potential as
  \begin{align}
      I_1(z) = \sum_{k=-K}^K c_k \oint_\Gamma \frac{\tau^k \dif \tau}{\tau - z},
  \end{align}
  where the integrals can be exactly evaluated using residue calculus,
  just as is done above for the unit circle. This will in fact work
  perfectly fine for some geometries (the unit circle in particular),
  but the interpolation problem can be severely ill-conditioned in a
  way that is hard to control, since it depends on the geometry
  itself. This could potentially be investigated further, but that is
  deemed beyond the scope of this work.
\end{rmk}

\subsection{Higher order integrals}

For higher order integrals ($m > 1$), we have
\begin{align}
  I_m(z) &= \oint_\Gamma \frac{\sigma \dif \tau}{\pars{\tau - z}^m}
           = \int_0^{2\pi} \frac{\sigma(t)\gamma'(t)\dif t}{(\gamma(t)-z)^m}
           ,
           \quad m \in \mathbb Z^+.
\end{align}
Evaluating this integral with SSQ is completely analogous to the case
when $m=1$, with the difference that we instead regularize with
$(e^{it}-e^{i\preimage})^m$, such that
\begin{align}
  f(t, \preimage) = \sigma(t)\gamma'(t)
  \pars{\frac{e^{it}-e^{i\preimage}}
  {\gamma(t)-\gamma(\preimage)}}^m
  .
\end{align}
We then need expressions for the integrals
\begin{align}
  p_k^m(\preimage) = 
  \frac{1}{i}
  \oint_S \frac{\xi^{k-1} \dif\xi}{\pars{\xi-\zeta}^m} .
\end{align}
Evaluating the residues in the same way as for the $m=1$ case, we get
the general expression for integer $m \ge 1$,
\begin{align}
  p_k(\preimage) = 
  2\pi\frac{\prod_{j=1}^{m-1} (k-j)}{(m-1)!} 
  e^{i(k-m)\preimage}
  \begin{cases}
     1, & \text{if } \Im \preimage>0 \text{ and } k \ge m,\\
     -1, & \text{if } \Im \preimage<0 \text{ and } k \le 0,\\
    0, & \text{otherwise}.
  \end{cases}
\end{align}

\subsection{Log kernel}

Consider now the log integral, also known as the Laplace single-layer
potential, for a real-valued density $\sigma$,
\begin{align}
  I_L(z) &= \oint_\Gamma \sigma \log|\tau-z|\:|\dif \tau| \\
  &= \int_0^{2\pi} 
    \underbrace{\sigma(t) |\gamma'|}_{f(t)}
    \log|\gamma(t)-z|\dif t,
\end{align}
where $f$ is assumed to be a smooth function. Proceeding in a similar
fashion as before, we find $\preimage$ and separate the integral as
\begin{align}
  I_L(z) = \int_0^{2\pi} f(t) \pars{
  \log|\gamma(t)-z|
  -
  \log|e^{it}-e^{i\preimage}|
  }\dif t
  + \int_0^{2\pi} f(t)\log|e^{it}-e^{i\preimage}| \dif t .
\end{align}
The first integral is now regularized, and can be evaluated using the
trapezoidal rule. The second integral we rewrite using a truncated
Fourier series and the relation $\log|r|=\Re\log r$, such that we can
evaluate it using SSQ,
\begin{align}
  \int_0^{2\pi} f(t)\log|e^{it}-e^{i\preimage}| \dif t \approx 
  \Re \sum_{k=-N/2}^{N/2} \hat f_k
  \underbrace{
  \int_0^{2\pi} e^{ikt}\log(e^{ikt}-e^{ik\preimage})\dif t
  }_{q_k(\preimage)} =
  \bm{\hat f} \cdot \bm q.
\end{align}
Just as with $p_k$, the integrals $q_k$ can be transformed into
integrals on the unit circle,
\begin{align}
  q_k(\preimage) = \frac{1}{i} \int_S \xi^{k-1}\log(\xi-\zeta)\dif\xi .
  \label{eq:qint}
\end{align}
Note that this is not necessarily a closed contour integral, as there
is a branch cut in the complex logarithm that needs to be taken into
account. We will now show to evaluate $q_k(\preimage)$ depending on the sign
of $\Im \preimage$.

\subsubsection{$\Im \preimage<0$}
Beginning with the case $|\zeta|>1$ (or $\Im \preimage<0$), we note that
$\xi-\zeta \ne 0$ on the entire unit disc. We can therefore choose the
branch cut of the complex logarithm such that it never intersects the
unit circle. This choice allows us to evaluate the integral in
\eqref{eq:qint} as a contour integral on the unit circle, with a
possible pole of order $1-k$ at the origin,
\begin{align}
  \frac{1}{2\pi i} \oint_S \xi^{k-1}\log(\xi-\zeta)\dif\xi 
  =
  \resbrac{\xi^{k-1}\log(\xi-\zeta), 0} 
  =
    \begin{cases}
      \zeta^{k} / k & \text{ if } k < 0, \\
      \log(-\zeta) & \text{ if } k=0, \\
      0 & \text{ if } k>0.
    \end{cases}          
          \label{eq:logcont_imneg}
\end{align}
Since $\hat f_0 \in \mathbb R$, we only need the real part of $q_k$ at
$k=0$.
We then have that  $\bm q$ is for $\Im \preimage<0$ computed as
\begin{align}
  q_k(\preimage) =
  2\pi
  \begin{cases}
    \frac{1}{k} e^{i k \preimage} & \text{ if } k<0, \\
    -\Im \preimage & \text{ if } k=0, \\
    0 & \text{ if } k>0 .
  \end{cases}
        \label{eq:qk_imneg}
\end{align}

\subsubsection{$\Im \preimage>0$}
We are now left with the case $|\zeta|<1$ (or $\Im \preimage>0$). Here,
$\xi-\zeta=0$ at at point on the unit disk, so the integral of
$\log(\xi-\zeta)$ on the unit circle must intersect the branch cut in
the $\xi$ plane going from $\zeta$ to infinity. We choose the branch
cut such that it passes through the point $\xi=1$. We define $1^+$ and
$1^-$ to be the limits approaching $1$ from above and below in the
plane, such that
\begin{align}
  \log(1^--\zeta)=\log(1^+-\zeta)+2\pi i.
  \label{eq:bcut_diff}
\end{align}
Then the integral on the unit circle $S$ is,
\begin{align}
  q_k(\preimage) = \frac{1}{i} \int_{1^+}^{1^-} \xi^{k-1}\log(\xi-\zeta)\dif\xi .
\end{align}
In order to evaluate this, we let $C$ be the closed contour following
the unit circle from $1^+$ to $1^-$, then going from $1^-$ to $\zeta$
on the negative side of the branch cut, and then going from $\zeta$ to
$1^+$ on the positive side of the branch cut,
\begin{align}
  \oint_C = \int_{1^+}^{1^-} + \int_{1^-}^{\zeta} + \int_{\zeta}^{1^+} .
\end{align}
The integrand is analytic inside $C$, except for a pole of order $1-k$
at the origin when $k \le 0$, with residue given by
\eqref{eq:logcont_imneg}. The integrals along the branch cut are
\begin{align}
  \int_{1^-}^{\zeta} + \int_{\zeta}^{1^+}
  &=
  \int_{1}^{\zeta} \xi^{k-1}
  \pars{\log(\xi-\zeta) + 2\pi i}\dif\xi + 
  \int_{\zeta}^{1} \xi^{k-1}\log(\xi-\zeta)\dif\xi \\
  &=
    2\pi i
    \begin{cases}
      \log(\zeta) & \text{ if } k=0,\\
      (\zeta^k - 1) / k & \text{ otherwise}.
    \end{cases}
\end{align}
Combining the results, and just as before keeping only the real part
of the $k=0$ integral, we arrive at the final results for $\Im \preimage>0$,
\begin{align}
  q_k(\preimage) =
  \frac{2\pi}{k}
  \begin{cases}
    1 & \text{ if } k<0, \\
    0 & \text{ if } k=0, \\
    (1 - e^{i k \preimage}) & \text{ if } k>0 .
  \end{cases}
        \label{eq:qk_impos}
\end{align}

\section{Experiments}
\label{sec:experiments}

In order to demonstrate the method outlined above, we here report
experimental results for the by now well-known ``starfish'' geometry
$\gamma(t) = (1 + 0.3 \cos 5t) e^{it}$. The method has been
implemented in Julia, and is available online as open source together
with the scripts that produce the numerical in this
manuscript\footnote{\url{https://github.com/ludvigak/TrapzSSQ.jl}}.

\subsection{Laplace double layer potential}
\label{sec:laplace}

As a first demonstration, we reuse the Laplace test problem from
\cite{AfKlinteberg2020line}. That is, we evaluate the solution to the
Laplace equation $\Delta u=0$ in the domain bounded by $\Gamma$, with
boundary condition $u=u_e(z)=\log\abs{3+3i-z}$. The solution to this
problem can be represented using the Laplace double layer potential (DLP)
$u(z)=\Im I_1(z)$, which leads to a second-kind integral equation in
the real-valued density $\sigma$. Once this integral equation is
solved, the layer potential $u(z)$ can be evaluated in the domain
using quadrature, and compared to the exact solution $u_e(z)$.

\begin{figure}[h]
  \centering
  \begin{subfigure}{.49\textwidth}
    \includegraphics[width=\textwidth]{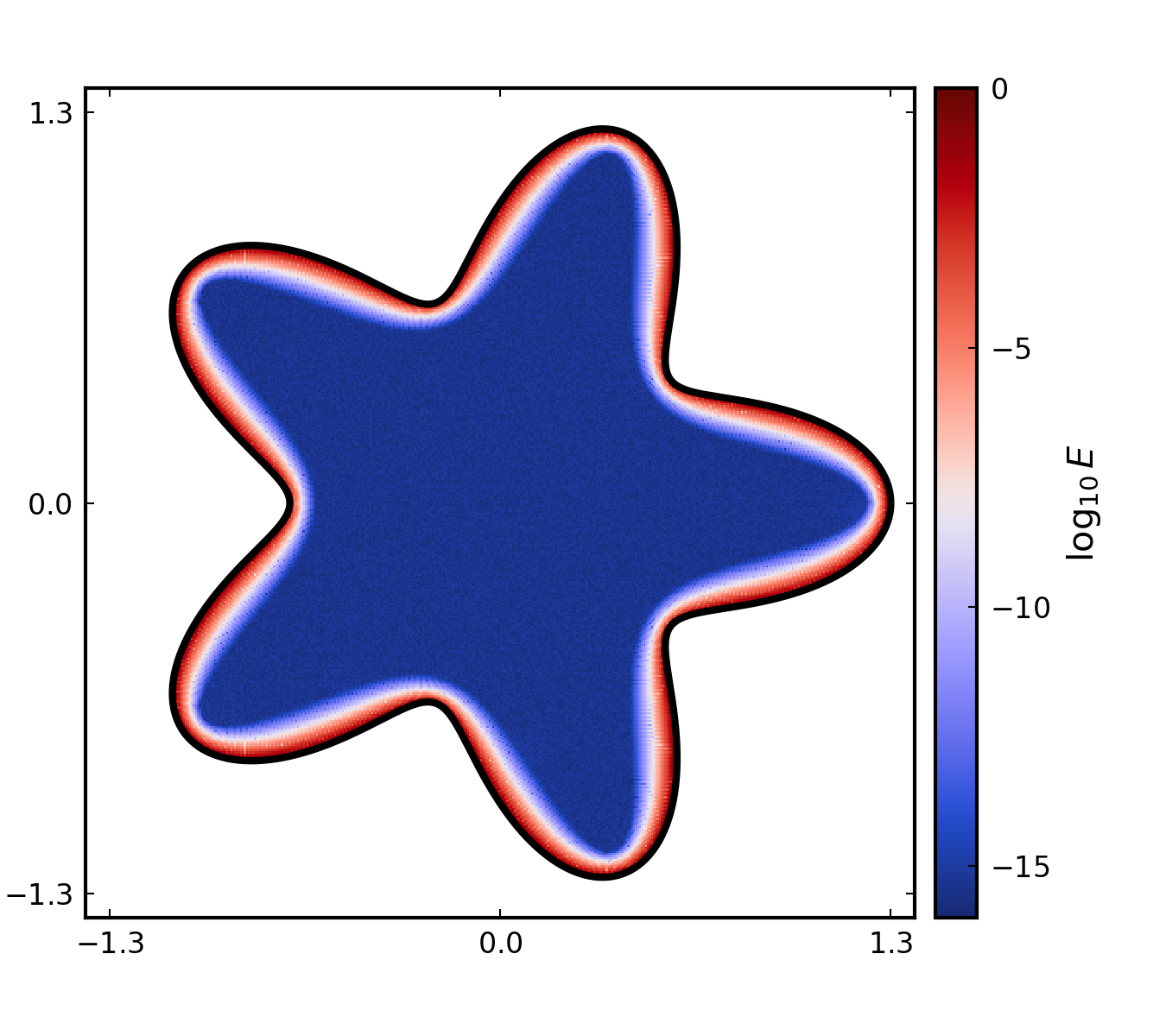}
    \caption{Trapezoidal rule, $E_{\max}=\mathcal O(1)$}
    \label{fig:laplace1}
  \end{subfigure}
  \hfill
  \begin{subfigure}{.49\textwidth}
    \includegraphics[width=\textwidth]{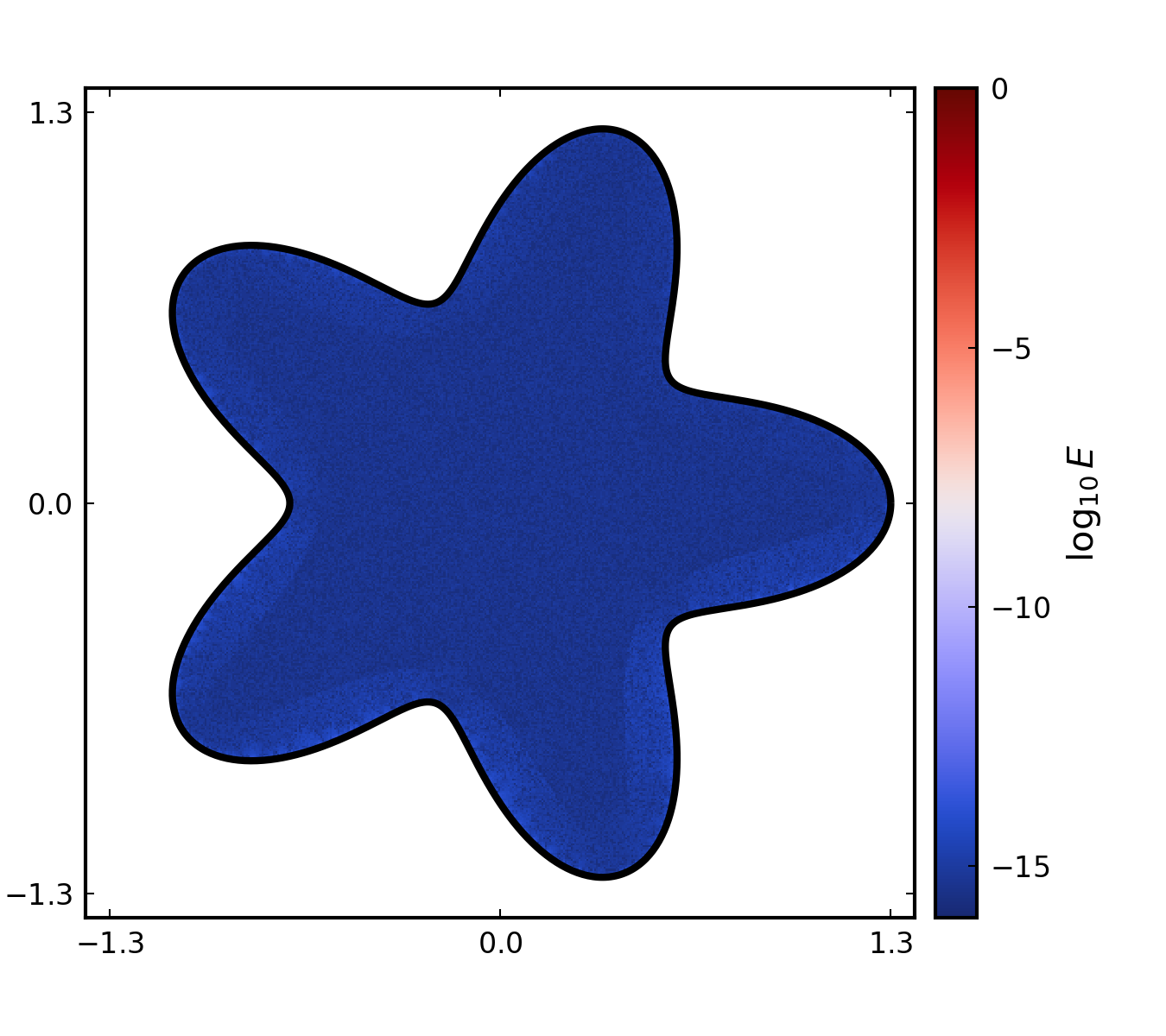}
    \caption{SSQ, $E_{\max}=7\cdot10^{-13}$}
    \label{fig:laplace2}
  \end{subfigure}  
  \caption{Error in the Laplace double layer potential on a domain
    with $N=400$ points on the boundary, comparison between the
    trapezoidal rule and SSQ.}
  \label{fig:laplace}
\end{figure}

In \cref{fig:laplace} we show results for the Laplace DLP when
$\Gamma$ is discretized using $N=400$ points that are equidistant in
the parametrization $t$. The layer potential is evaluated using the
trapezoidal rule on a $400 \times 400$ grid inside the domain, but for
points close to the boundary we also evaluate using SSQ and compare
the results. As expected, evaluation using the trapezoidal rule leads
to $\mathcal O(1)$ errors near the boundary, while SSQ is accurate at
all points where it is used. However, full machine precision accuracy
is not recovered; the maximum error in SSQ close to the boundary is
$\mathcal O(10^{-13})$.

\subsection{Convergence}
\label{sec:convergence}

For a more systematic investigation of the method, we create the
following setup. Let $\Gamma$ be the starfish domain from the previous
section, and discretize it using $N$ points. We then create 100
evaluation points as $z=\gamma(\preimage)$, where
$\Re\preimage \in [0, 2\pi)$ and $\Im\preimage=d$. In order to test
points both interior and exterior to $\Gamma$, we use
$d = \{\pm 0.01, \pm 0.02, \pm 0.04\}$. For a range of $N$, we then
compute the integrals $I_L, I_1, I_2, I_3$ using both trapezoidal and
singularity swap quadrature, and for each $N$ save the maximum error
across all $z$. The density $\sigma$ and reference solution used are:

\begin{figure}[p] \centering
  \begin{subfigure}{.49\textwidth}
    \includegraphics[width=\textwidth]{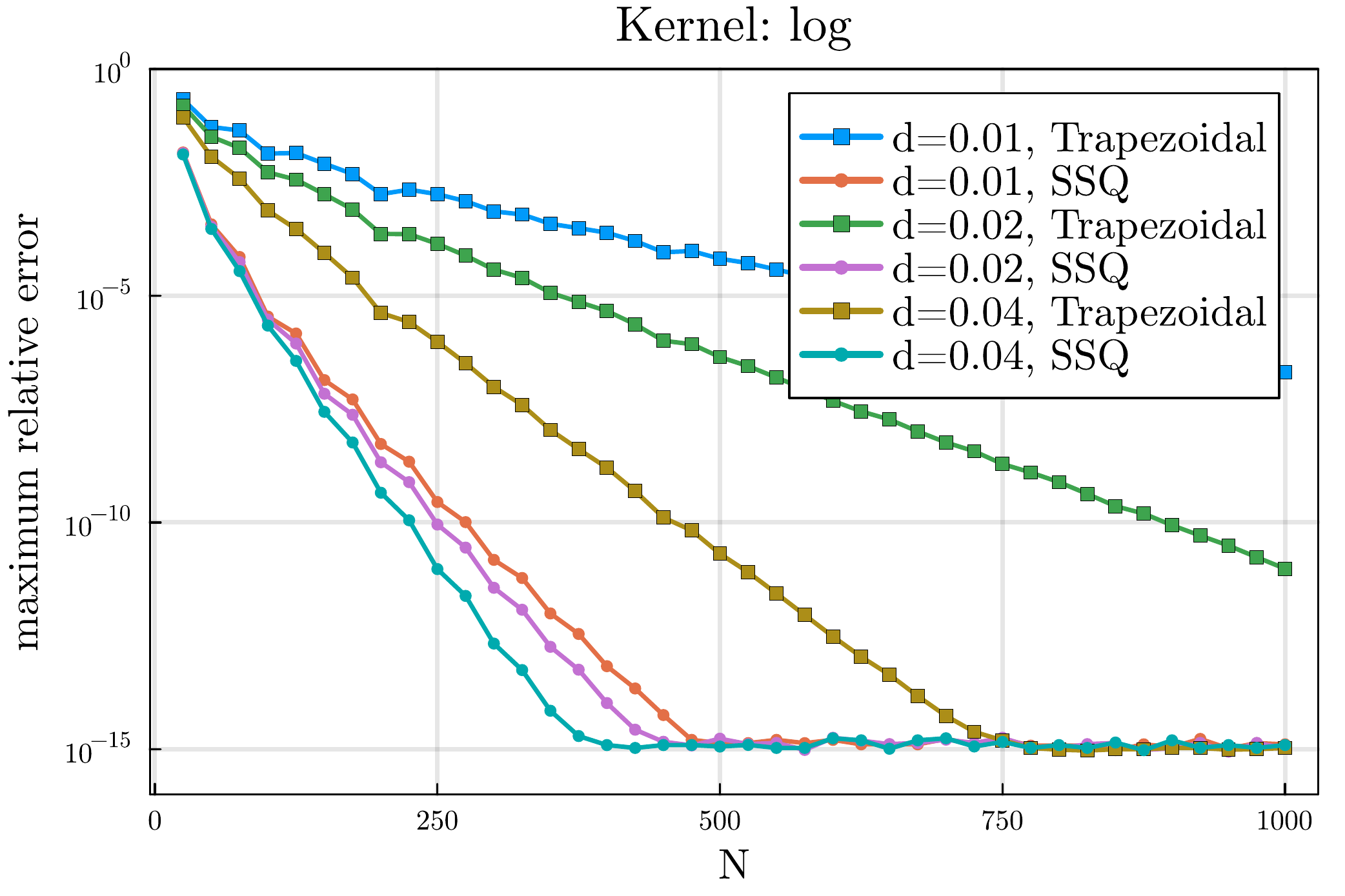}
    \caption{}
    \label{fig:convlog_interior}
  \end{subfigure}
  \hfill
  \begin{subfigure}{.49\textwidth}
    \includegraphics[width=\textwidth]{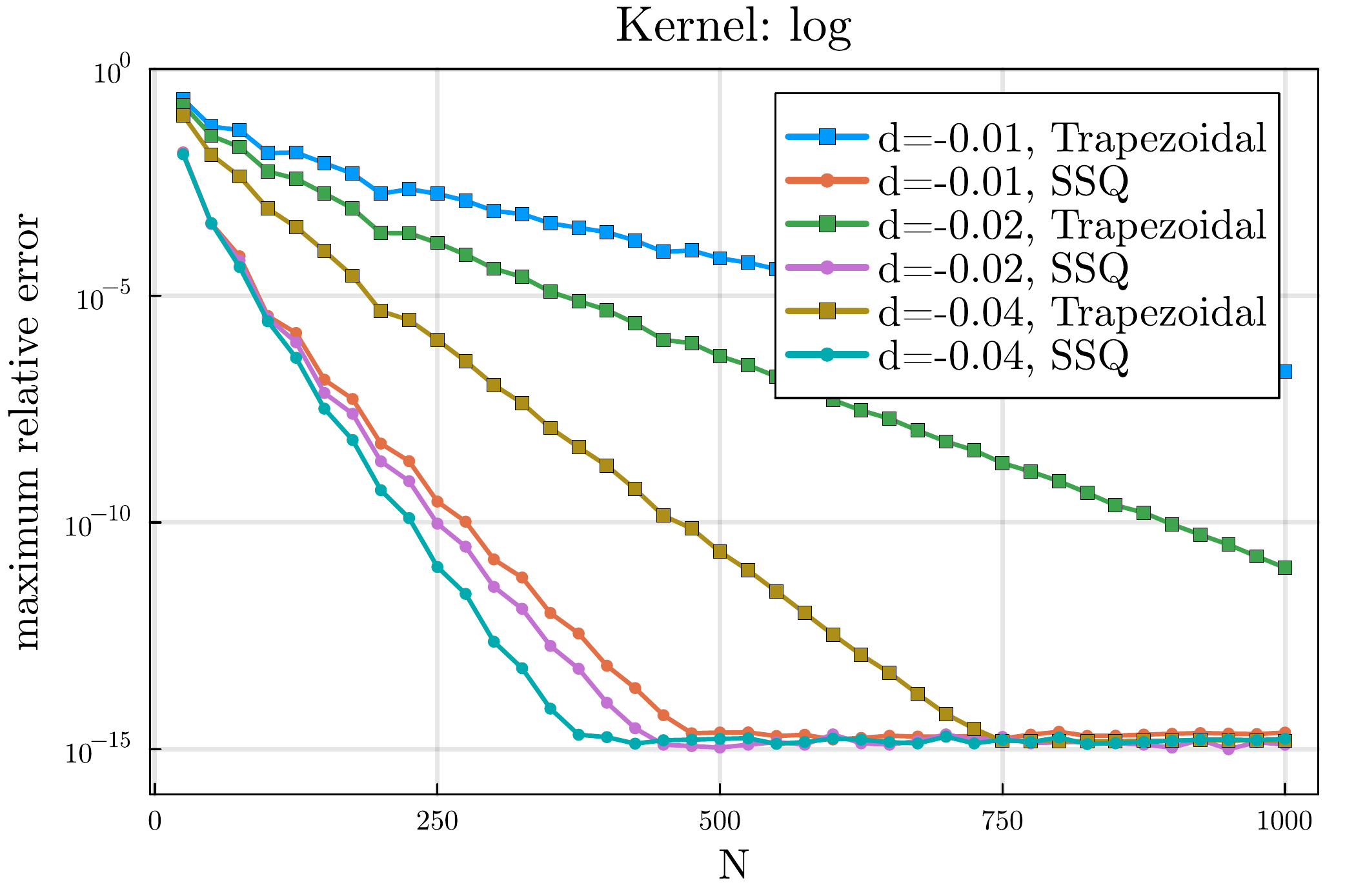}
    \caption{}
    \label{fig:convlog_exterior}
  \end{subfigure}
  \\
  \begin{subfigure}{.49\textwidth}
    \includegraphics[width=\textwidth]{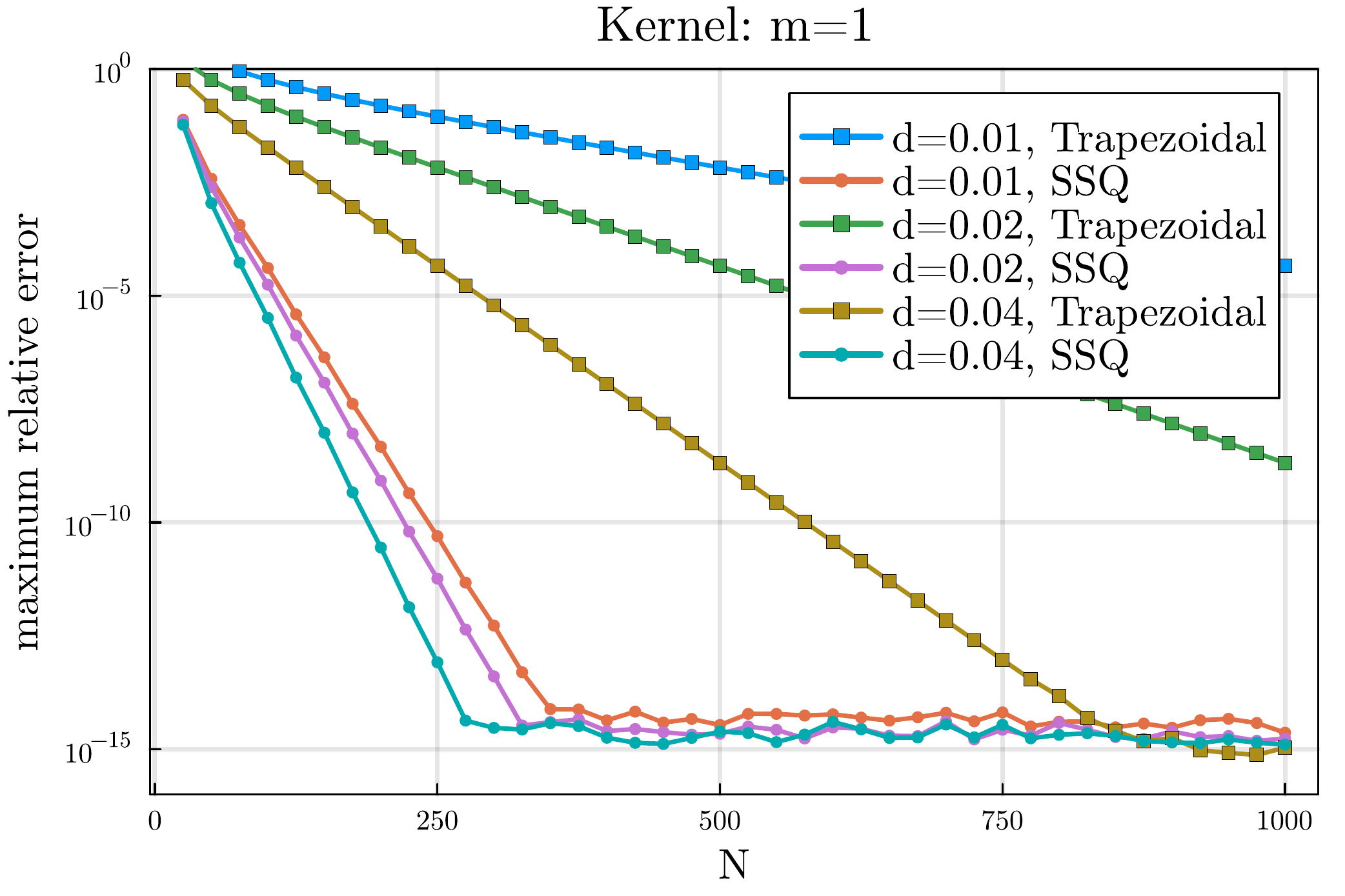}
    \caption{}
    \label{fig:conv1_interior}
  \end{subfigure}
  \hfill
  \begin{subfigure}{.49\textwidth}
    \includegraphics[width=\textwidth]{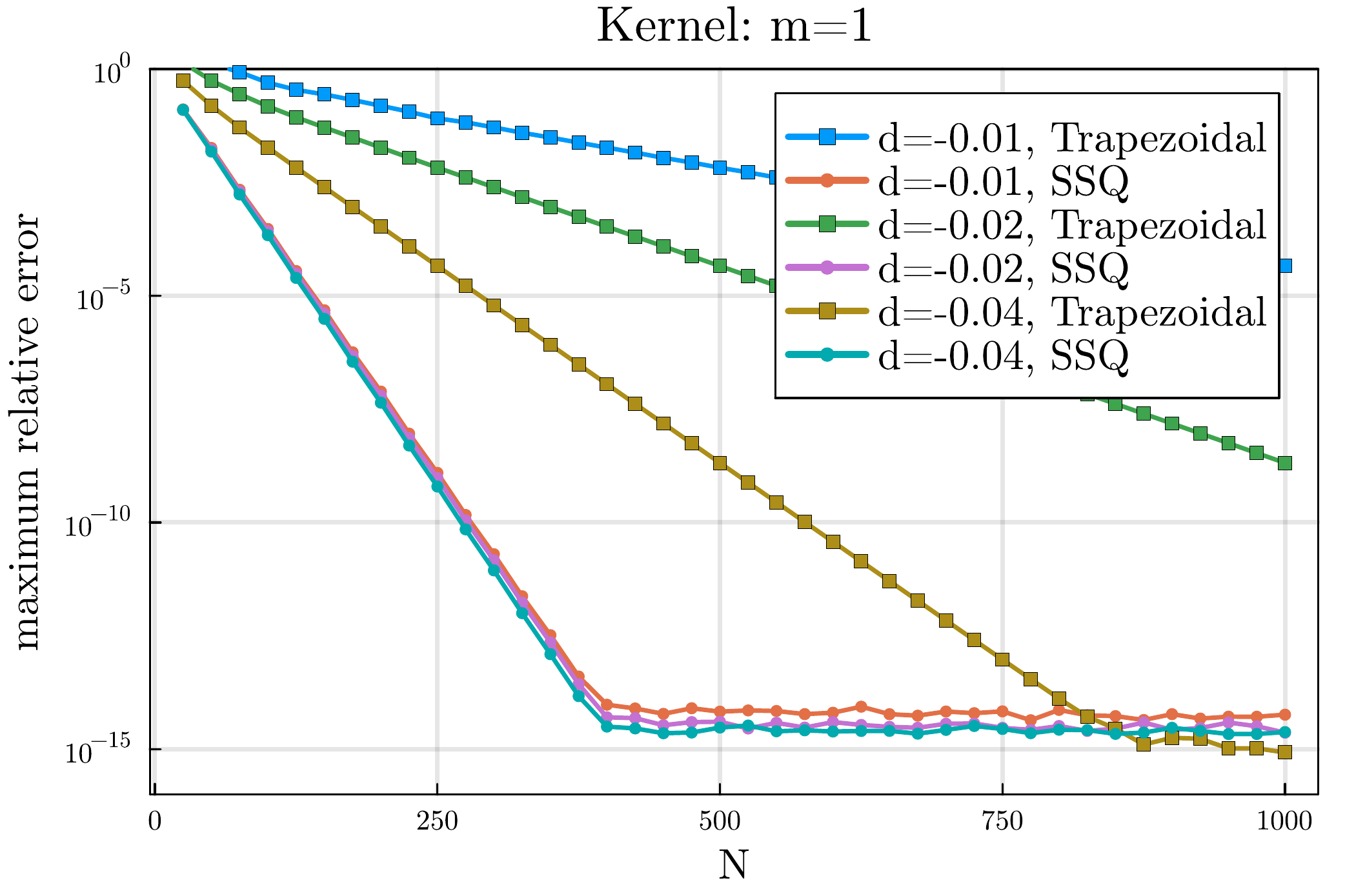}
    \caption{}
    \label{fig:conv1_exterior}
  \end{subfigure}
  \\
  \begin{subfigure}{.49\textwidth}
    \includegraphics[width=\textwidth]{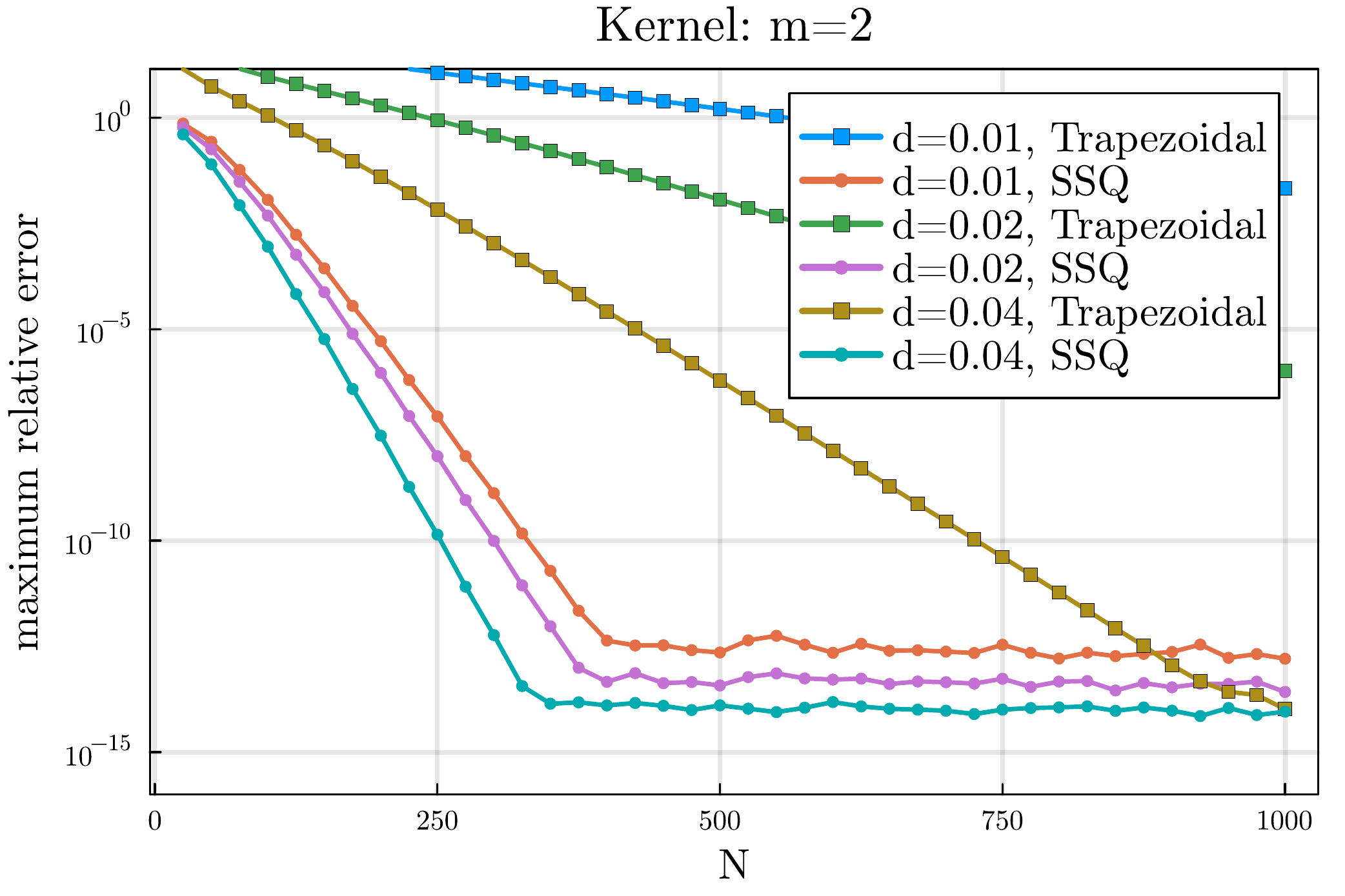}
    \caption{}
    \label{fig:conv2_interior}
  \end{subfigure}
  \hfill
  \begin{subfigure}{.49\textwidth}
    \includegraphics[width=\textwidth]{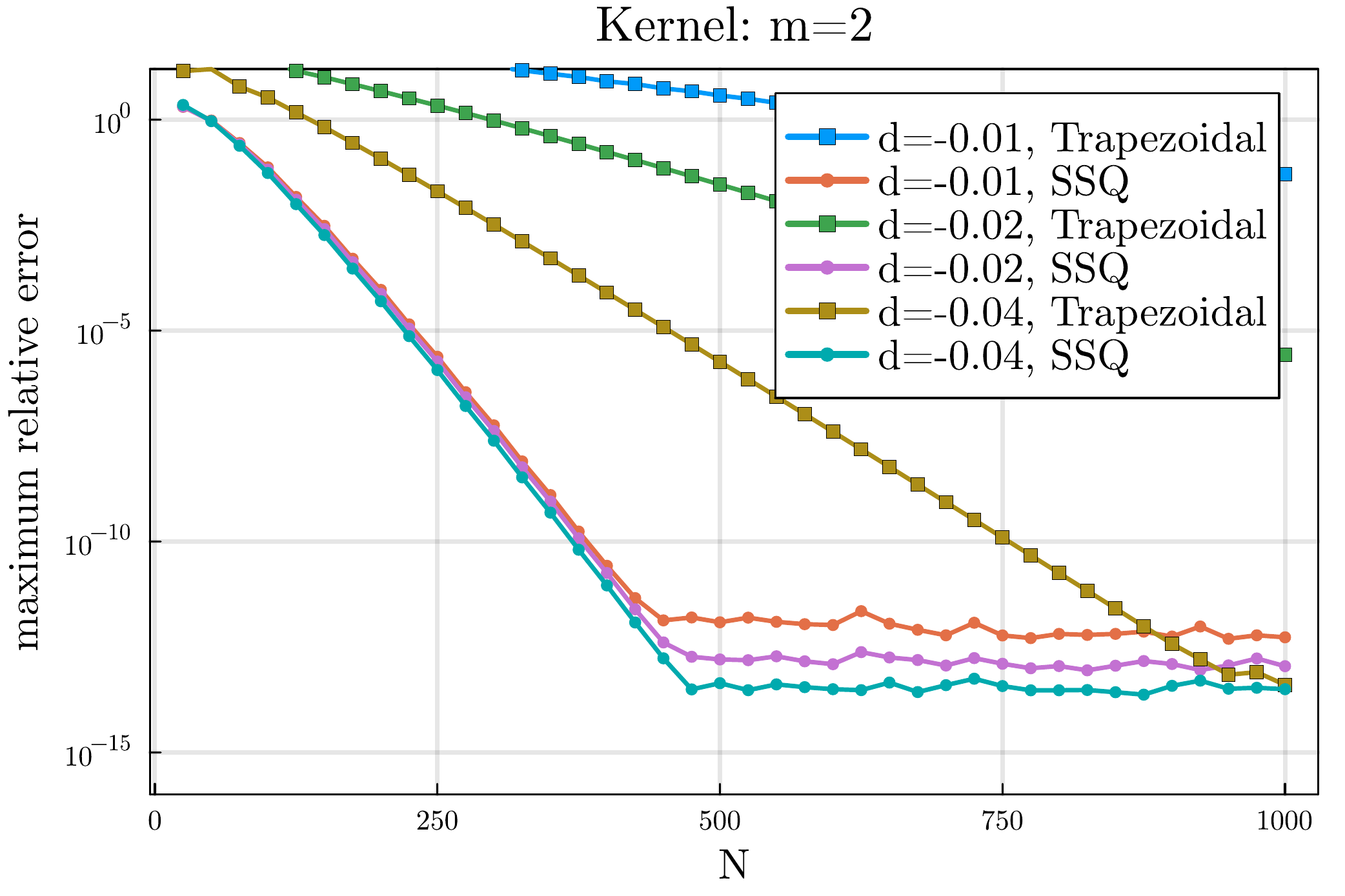}
    \caption{}
    \label{fig:conv2_exterior}
  \end{subfigure}
  \\
  \begin{subfigure}{.49\textwidth}
    \includegraphics[width=\textwidth]{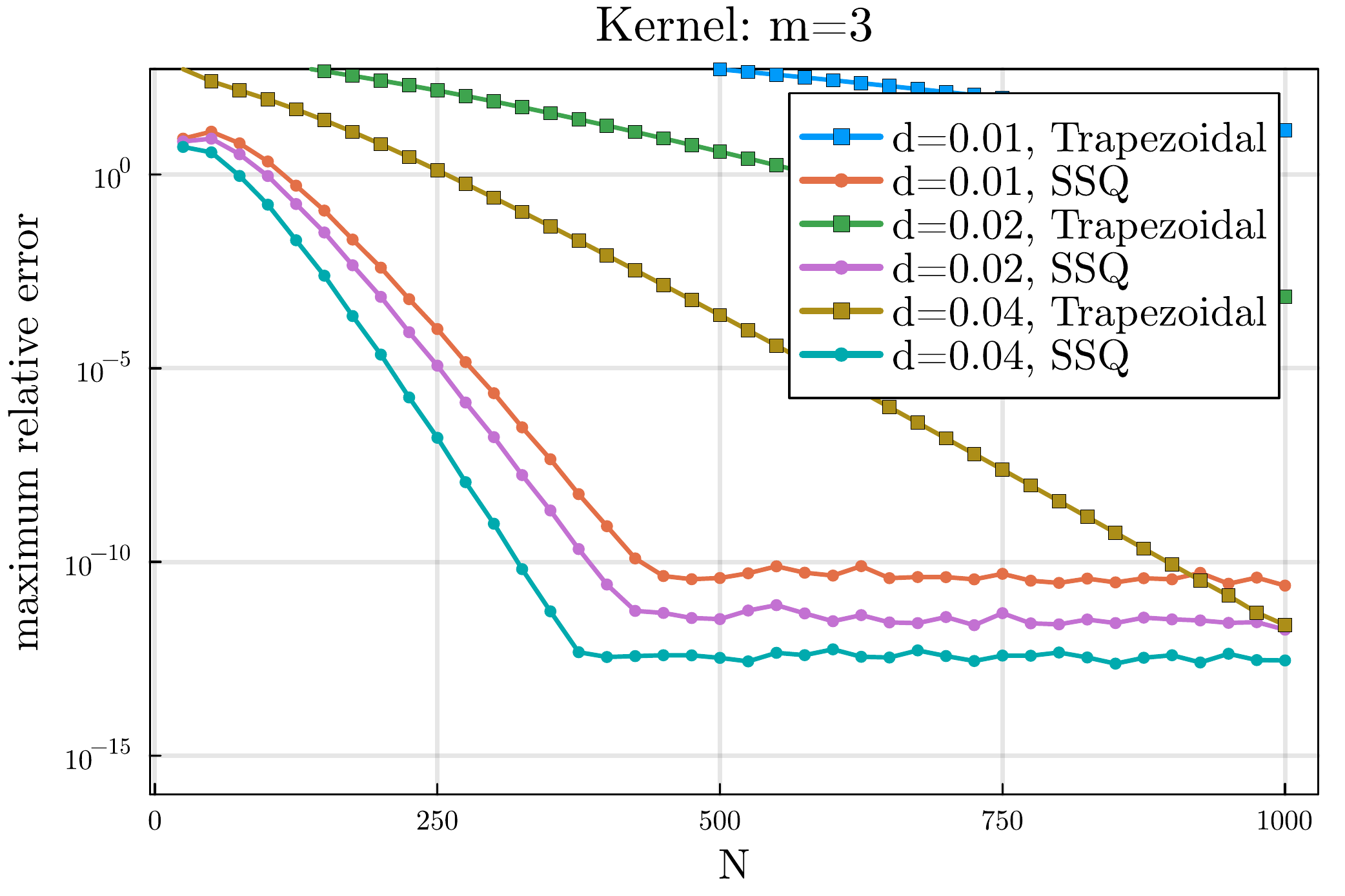}
    \caption{}
    \label{fig:conv3_interior}
  \end{subfigure}
  \hfill
  \begin{subfigure}{.49\textwidth}
    \includegraphics[width=\textwidth]{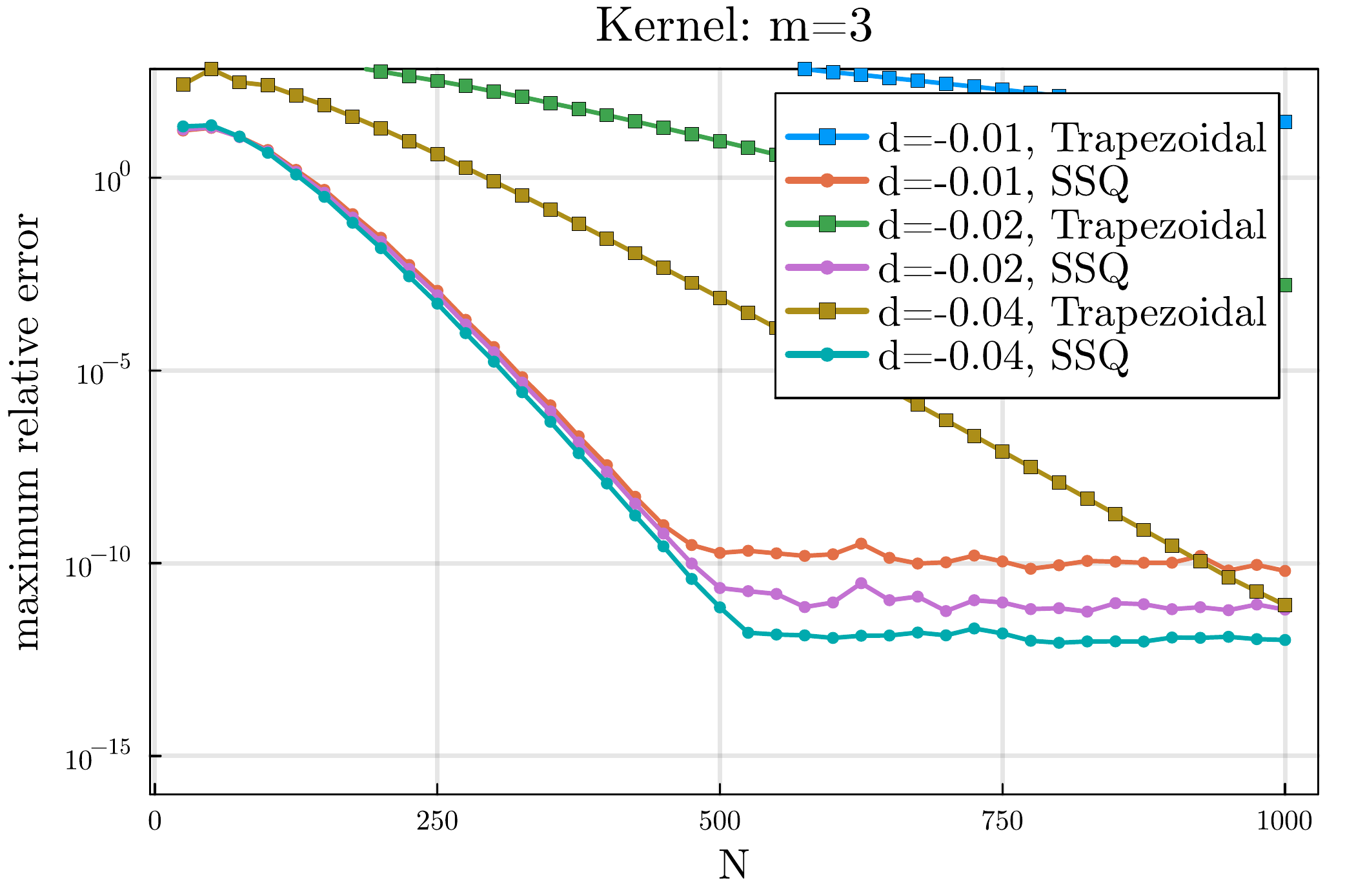}
    \caption{}
    \label{fig:conv3_exterior}
  \end{subfigure}
  \caption{Convergence of the trapezoidal rule and SSQ for the log
    kernel and $m=1,2,3$, measured as the maximum error on a set of
    target points such that $\Im\preimage=d$. Left column shows
    interior points, right column shows exterior points.}
  \label{fig:conv}
\end{figure}

\begin{itemize}
\item For $I_L$ we use the $\sigma(t) = g_1(t)g_2(t)$ and compute the
  reference solution using adaptive Gauss-Kronrod quadrature.
\item For $I_m$ and \emph{interior} points ($\Im\preimage>0$) we use
  $\sigma(\tau)=\tau^3 + \tau$, with reference solution $I_m(z)=2 \pi i \sigma^{(m)}(z)/(m-1)!$.
\item For $I_m$ and \emph{exterior} points ($\Im\preimage<0$) we use
  $\sigma(\tau)=\tau^{-1}$, with reference solution $I_m(z)=2\pi i z^{-m}$.
\end{itemize}

The results of this investigation are shown in \cref{fig:conv}. As
expected, the convergence rate of the trapezoidal rule depends
strongly on the distance to the boundary, which is not the case for
SSQ. There appears to be a weak dependence on the rate for interior
points, but that could be an effect of the chosen test problem. More
concerning is that lowest attainable error appears to increase for
smaller $d$ when $m>1$. It is unclear what causes this effect,
although in experiments there appears to be an upper bound to the
error growth.

\section{Conclusions}

We have extended the singularity swap quadrature (SSQ,
\cite{AfKlinteberg2020line}) method to closed 2D curves discretized
using the trapezoidal rule. This extension builds on the simple
observation that interpolatory quadrature can be used on a periodic
integrand if the problem is first swapped to the unit circle. The
method relies on representing quantities as Fourier expansions, and as
a consquence the domintaing cost of the method is that of applying the
fast Fourier transform (FFT).

The method is accurate, exponentially convergent, and relatively
simple to implement (and we provide source code). The fact that the
cost is $\mathcal O(N \log N)$ per target point makes the method
unfeasible for large problems, but that is a consequence of the
underlying trapezoidal rule being global (and exponentially
convergent). For maximum efficiency, composite Gauss-Legendre
quadrature with the original SSQ likely remains the better option.

In the original SSQ paper, it was shown that the method can be applied
to line integrals in 3D by finding the singularity preimage through
analytic continuation of the squared-distance function. The same
technique can be applied to the present method, for closed curves in
three dimensions. What remains is analytical evaluation of the
integrals of the interpolatory quadrature. This is straightforward for
the logarithmic kernel, but the integrals corresponding to power-law
kernels are more challenging. Deriving expressions for these is
currently work in progress.

\subsection*{Acknowledgements}
  The key contributions of this work were made with support from the
  Knut and Alice Wallenberg Foundation, under grant no.\ 2016.0410.

\clearpage

\bibliographystyle{abbrvnat_mod}
\bibliography{library}

\end{document}